\documentstyle[mydefs,12pt,psfig]{article}

\newcommand{\be}{\begin{equation}}
\newcommand{\ee}{\end{equation}}
\newcommand{\la}{\left|\begin{array}{ccc}}
\newcommand{\ra}{\end{array}\right|}

\newcommand{\BE}{Biedenharn-Elliott relations}

\newcommand{\eins}{{\rm \setlength{\unitlength}{1em}
\begin{picture}(0.75,1)
\put(0,0){1}\put(0.34,0){\line(0,1){0.65}}
\end{picture}}}
\begin{document}

\input amssym.def
\input amssym.tex

\begin{center}
{ \LARGE \bf On the relation between two quantum group invariants
of 3-cobordisms
%\\[0.25cm]
%topological quantum field theories
% \\and invariants of graphs\\
}
 \\[1.5cm]

\begin{tabular}{cc}
Anna Beliakova
%\footnotemark[1]
& Bergfinnur Durhuus\\
 Institut f\"{u}r Theoretische Physik &   Mathematics Institute\\
 Freie Universit\"{a}t Berlin
& University of Copenhagen\\
\end{tabular}
\vspace*{0.5cm}
\\January 95 \vspace*{0.5cm}

\end{center}

\parbox{10cm}{{\footnotesize {\bf Abstract}:  We prove in the context
of quantum groups at even roots of unity that a Turaev-Viro type
invariant of a 3-dimensional cobordism $M$
equals the tensor product of the Reshetikhin-Turaev invariants of $M$
and $M^\ast$, where the latter denotes M with orientation reversed.
}}

\section{Introduction}\label{ein}
%\footnotetext[1]{Supported by DAAD and DFG, SFB 288
%''Differentialgeometrie und Quantenphysik``}
%A topological quantum field theory in dimension 2 (TQFT) [At]
%provides a
%bridge between topology, link invariants  and quantization
%of Chern-Simons theory.
%More exactly,TQFT
%associates a
According to [At] a 3-dimensional topological quantum field
theory (TQFT) associates a finite dimensional vector space $V_\Sigma$ to each
compact closed oriented
2-dimensional surface $\Sigma$ and a vector (partition function) $Z(M)\in
V_\Sigma$ to each
compact oriented 3-dimensional manifold $M$ with boundary
$\Sigma$, satisfying a certain set of axioms. Of particular relevance
for the following discussion are the following: 1) $V_{\Sigma^\ast}$ is
the dual space of $V_\Sigma$ for each surface $\Sigma$, where
$\Sigma^\ast$ denotes $\Sigma$ with orientation reversed, 2) given an
orientation preserving diffeomorphism $f: \Sigma \rightarrow
\Sigma^\prime$  between oriented surfaces, there exists an isomorphism
$U(f): V_{\Sigma} \rightarrow V_{\Sigma^\prime}$ fulfilling $U(f_1f_2)
= U(f_1)U(f_2)$ for any pair of diffeomorphisms that can be composed, and  3)
if M
is obtained by gluing two 3-manifolds $M_1$ and $M_2$ along
$\Sigma\in \partial {M_1}$ and $\Sigma^\ast\in \partial {M_2}$ then
$Z(M)$ is obtained by contracting $ Z(M_1)\otimes Z(M_2)$ with respect
to $V_\Sigma$. In addition, the vectorspace associated to the empty
surface is assumed to be the complex numbers, and if $\Sigma$ is the
disjoint union of two surfaces $\Sigma_1$ and $\Sigma_2$ then
$V_\Sigma = V_{\Sigma_1} \otimes V_{\Sigma_2}$. In particular, if $M$ is a
closed manifold  $Z(M)$ is a complex number which is a topological
invariant of $M$.
%In case $M$ has a non-empty boundary we may more generally view $Z(M)$
%as a cobordism invariant.

Alternatively, the gluing property 3) can be reformulated in
terms of operators as follows. Viewing $M_1$ and $M_2$ as cobordisms
with $\partial{M_1} = \Sigma_1 \cup \Sigma$ and $\partial{M_2} =
{\Sigma^\prime}^\ast \cup \Sigma_2$ we can correspondingly consider
the state sums as operators $Z(M_1): {V_{\Sigma_1}}^\ast
 \rightarrow
V_\Sigma$ and $Z(M_2): V_{\Sigma^\prime} \rightarrow V_{\Sigma_2}$ by 1).
Given an orientation preserving diffeomorphism $f: \Sigma \rightarrow
\Sigma^\prime$ and letting $M$ denote the manifold obtained by gluing
$M_1$ to $M_2$ along $f$,  property 3) is equivalent to
\be \label{glu} Z(M) = Z(M_2) U(f) Z(M_1)\;  .\ee
Note that the symmetry of the gluing w.r.t. $M_1$ and $M_2$ requires
that
\be \label{transu}  U(f^\ast) = (U(f)^t)^{-1}, \ee
 where $f^\ast: \Sigma^\ast
\rightarrow (\Sigma^\prime)^\ast$ denotes $f$ with orientations on
$\Sigma$ and $\Sigma^\prime$ switched, and the superscript $t$
indicates transposition.
There now exists in the literature a variety of rigorous constructions of
3-dimensional TQFT's. In this note we shall consider
the constructions by Reshetikhin-Turaev [RT] and the one by
Turaev-Viro [TV] and their generalizations (see [T], [DJN], [KS],
[BD]). These are all based on the algebraic structure of the
representation theory of quantum groups with deformation parameter
equal to a root of unity, and are known to be related to
Chern-Simons theory with an arbitrary compact gauge group.

In [BD] we have proven that for closed manifolds the invariant
$Z_{TV}(M)$ of the Turaev-Viro construction equals the modulus squared
of the invariant $\tau(M)$ obtained by the Reshetikhin-Turaev construction
for a general quantum group at simple even roots of unity (see also
[Wa], [T] and [R]). The purpose of this paper is to extend this result to
manifolds with boundary, i.e. we show that
$$Z_{TV}(M)= \tau(M)\otimes \tau(M^\ast)$$
for any 3-cobordism $M$. Here $Z_{TV}(M)$ and $\tau(M)$ denote the
cobordism invariants defined in [BD] and [T],
respectively. In section 2 we recall briefly the basic elements of the
Turaev-Viro construction as developed in [BD] and refer the reader to
that paper for fuller details. We then prove a basic lemma which
yields certain isomorphisms from the state spaces of the theory
onto certain explicitly realizable spaces. This result is used in
Section 3 to obtain an equivalent TQFT for which the announced
factorization property is then proven.

\newtheorem{satz}{Theorem}
\newtheorem{lem}[satz]{Lemma}

\section {Turaev-Viro TQFT}\label{sixj}
In this section we briefly recall the formulation and basic properties of TQFT
of the Turaev-Viro type (for more details see [BD]).
The corresponding  state sum will be denoted
by $Z(M)$ (omiting the index $TV$ in the following).

Originally, the Turaev-Viro invariant was defined for a compact
connected closed
oriented 3-manifold $M$ as follows [TV]:
Choose a triangulation of $M$ and associate to each 1-simplex of the
triangulation an index (or a colour) from a finite set $\cal{I}$
of so-called ``admissible'' representations of a quantum group.
To each coloured tetrahedron one then associates a 6j-symbol,
which is possible due to the invariance of 6j-symbols under the
tetrahedral symmetry group. In addition, to each coloured link one
attaches a factor $\omega^2_i$, which equals the quantum dimension of the
corresponding colour $i$, and to each vertex one attaches a factor
$\omega^{-2}$, where
$$\omega^2 = \sum_{i\in {\cal I}} \omega^4_i.$$ The invariant
$Z(M)$ is then obtained as the sum over all colourings of the
triangulation of  the product of all factors associated to tetrahedra,
links and vertices. It can be shown (using the \BE\ for 6j-symbols)
that the resulting quantity is independent of the particular choice of
triangulation.

We have here assumed that the 6j-symbols are scalars, i.e. that the
multiplicity of any representation $i\in{\cal I}$ in a tensor product
of two representations in $\cal I$ is always 0 or 1, which e.g. is the
case for $SU_q(2)$. For more general quantum groups the 6j-symbols are
tensors. To be specific we associate to
each oriented, coloured triangle $t$ in $ \Sigma =\partial M$ with
oriented boundary links as indicated in
Fig.1 (where the orientation of the plane is assumed to be
counter clock-wise)
 the vector space $V^k_{ij}$
of Clebsch-Gordan coefficients defined by $$H_i\otimes H_j =
\sum_{k\in {\cal I}}V^k_{ij}\otimes H_k,$$ where $H_i$ denotes the
vector space of the representation $i$.

\input epsf.sty
\begin{center}
\mbox{\epsfysize=2.5cm
\epsffile{dreinew1.ps}}
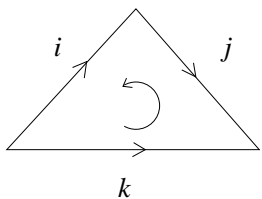
\\Fig.1 {\it An oriented $ \{i,j,k\}$-coloured  2-simplex}
\end{center}
%example, to the 2-simplex  on Fig.1
% corresponds the vector space $C^k_{ij}$
%$\{i,j,k\}\in {\cal I}^3$, where

The canonically dual vector space $(V^{k}_{ij})^\ast= V^{ij}_k$
will be associated to the oppositely oriented triangle.
For other configurations of arrows than that on Fig.1 the
corresponding spaces are defined by requiring that reversing an arrow
on a 1-simplex
is equivalent to
replacing its colour by the dual one (i.e. replacing the corresponding
representation by its adjoint).

%There is a canonical duality
%between the spaces associated to a
%given oriented coloured triangle and its oppositely oriented counter part.

Moreover, the 6j-symbol associated to an oriented
coloured tetrahedron with oriented links belongs to the tensor product
of the vector spaces associated to the triangles in its boundary.
Thus, we may define $Z(M)$ by replacing above the product of
6j-symbols by the corresponding tensor product and contracting with
respect to the dual pairs of spaces associated to triangles (with some
fixed orientation of links), and the result is again independent of the
choice of triangulation as well as of the chosen orientation of links.
In fact, this definition is easily extended to non-closed, oriented
manifolds $M$ by simply contracting only with respect to dual pairs of
spaces associated to $interior$ triangles of the triangulation. One then
obtains a tensor $Z'(M)$ in  the vector space $V'_{\partial
M}$ defined as the direct sum over all colourings of the links in
$\partial M$ of the tensor product of the spaces associated to the
triangles in $\partial M$. This space, of course, depends on the
triangulation of $\partial M$. However, any two such triangulations
may be connected by a triangulation of the cylinder $\partial M \times
[0,1]$ in the obvious sense, and $Z'(\partial M \times[0,1])$
defines a {\it cylinder map} between the corresponding spaces. In particular,
choosing the same triangulation at the two ends of the cylinder the
map becomes a projection, and the supports of the projections so
obtained may be canonically identified by the cylinder maps thus
defining the vector space $V_{\partial M}$, and at the same time
the partition functions $Z'(M)$ are also identified with a unique
vector $Z(M)\in V_{\partial M}$ fulfilling the required properties.

Exploiting ideas of Turaev [Tu] an effective calculational tool was
developed in [KS] by introducing coloured graphs $G_{\underline{x}}$
on the boundary of the manifold M and defining an associated
state sum $Z(M,G_{\underline{x}})$ generalizing $Z(M)$. Here a
coloured graph $G_{\underline{x}}$ is a closed 1-dimensional
simplicial complex, whose 0-simplexes have order at most 3 and whose
lines (i.e. maximal sequences of 1-simplexes joined by vertices of
order 2 are oriented and coloured (by elements in $\cal I$), the
collection of colours being indicated by $\underline{x}$. The graph
is assumed to be embedded into $\partial M$ such that over- and
undercrossings are distinguished.
The definition of $Z(M,G_{\underline{x}})$ proposed in [KS]
has the following
geometrical interpretation (see [BD]). One glues to the boundary
$\Sigma$ of $M$ a certain pseudo-manifold $P_G$ whose boundary
consists partly of one copy of $\Sigma^\ast$ (triangulated as
$\Sigma$) and partly of a surface on which the dual graph of $G$
determines a cell decomposition into triangles (corresponding to
3-vertices) and rectangles (corresponding to over- and
undercrossings) and whose edges inherit a colouring from
$\underline{x}$. The state sum $Z(M,G_{\underline{x}})$ then equals
$Z(M_{G_{\underline{x}}})$, where $M_{G_{\underline{x}}}$ is the
resulting pseudo-manifold with fixed colouring of boundary links given
by $\underline{x}$. Actually, the construction requires a slight
modification in case rectangles are present in the boundary (see
[BD]). Suffice here to mention that $Z(M,G_{\underline{x}})$ in all
cases belongs to the tensor product of the vector spaces associated to
the triangles dual to the 3-vertices in $G_{\underline{x}}$ and is a
homotopy invariant of the coloured graph $G_{\underline{x}}$ in
$\Sigma$.

In case $G$ is empty the pseudo-manifold $P_G$ is the cone over
$\Sigma$ and the boundary of the resulting manifold degenerates to a
point. On the other hand, if $G$ is sufficiently ``large'' so that $P_G$ is
homeomorphic to the cylinder $\Sigma \times [0,1]$, then
$M_{G_{\underline{x}}}$ is homeomorphic to $M$, and if $G$ in addition
has no over- or undercrossings it follows that
$\oplus_{\underline{x}}Z(M,G_{\underline{x}})$ equals $Z^\prime(M)$
with $\partial{M}$ triangulated by the dual graph to $G$.

The gluing axiom described at the beginning of section 1  can now be
reformulated in the language of graphs as follows. If
%According to the definition of the state sum
% gluing along $\Sigma$ can be described by the sum over all
%colourings of a canonical graph, i.e.\
$M$ is obtained by gluing  $M_1$ and $M_2$ along $\Sigma$ we have
\be Z(M)=\frac{1}{\omega^2}\sum_{\underline{x}}
\omega^2_{\underline{x}}\; Z(M_1, G^F_{\underline{x}})\; Z(M_2,
G_{\underline{x}}) \ee
%where  $\omega^2=\sum_{i\in {\cal
%I}}\omega^4_i$
%$\omega^2_i$ is the quantum dimension of the $i$-th
%irreducible representation.
for any canonical graph $G$ without over- or undercrossings, and where  $F$ is
 the gluing homeomorphism and $G^F$ denotes the image of $G$ under $F$.

The state sums $Z(M,G_{\underline{x}})$ satisfy a number of simple
relations under certain elementary changes of the graph
$G_{\underline{x}}$, which together with (2.1) can be used to show that
the dimension of $V_{\Sigma_g}$, where $\Sigma_g$ is a connected
surface of genus $g\geq 1$, is given by the square of the Verlinde formula:
%Our next task is to find the ``simpliest'' or minimal
%canonical graph on
%$\Sigma$,
%the number of whose colourings is equal exactly to the dimension of
%$V_\Sigma$. It would define a basis of the vector space $V_\Sigma$.
%The dimension of $V_\Sigma$ was calculated in [BD] (see also [KS])
%and is given by the square of the Verlinde formula [V]:
\be dim V_{\Sigma_g} = tr\; id_{V_{\Sigma_g}} =tr Z(\Sigma_g\times I) =
Z({\Sigma_g}\times S^1)=(tr {\vec{N}}^{2(g-1)})^2\ee
where ${\vec{N}}^2 = \sum_a (N^a)^2$ and $(N^a)$ is the multiplicity
matrix given by
\be (N^a)_{bc} = N^a_{bc} = dim V^a_{bc}\ee
for $a,b,c \in \cal I$.

It is even possible to realize the space $V_{\Sigma_g}$ explicitly as
follows. Consider a handle body $M_g$ of genus $g$ in $R^3$ with $\partial M_g
= \Sigma_g$ and introduce two copies  $c^L$ and $c^R$ of the graph
depicted below such that they are deformation retracts of $\Sigma$ in
$M_g$ and such that they are disjoint (and not linked).
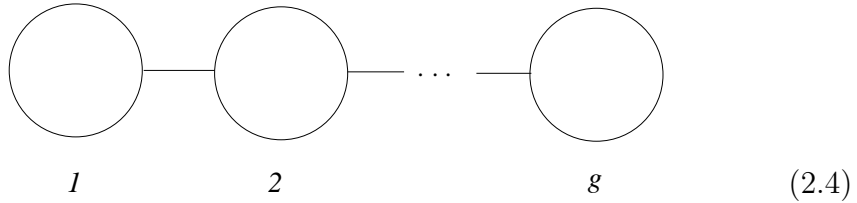
\be \label{kanon}
\mbox{\epsfysize=2.5cm\epsffile{kanonnew.ps}}
\ee
 Clearly
$c^L$ and $c^R$ then possess tubular neighborhoods that are disjoint
and diffeomorphic to $M_g$ and whose boundaries are
homotopic to $\Sigma_g$ in $M_g$. Removing two such tubular
neighborhoods from $M_g$ yields a manifold $\tilde{M_g}$ with three boundary
components $\Sigma_g$, $(\Sigma^{L}_g)^\ast$ and $(\Sigma^{R}_g)^\ast$
all of genus $g$.
We now choose the coordinates so that the cores $c^L$ and $c^R$ lie in the
$xy$-planes and their $z$-components  are equal  to {1} and
{-1} respectively. We will call the part
 of $\Sigma^L_g$ (resp. $\Sigma^R_g$)
where $z>1$ (resp. $z> -1$) the upper side and the other part
where $z<1$ (resp. $z< -1$)
the back side of $\Sigma^L_g$ (resp. $\Sigma^R_g$).

Next, we embed a copy $G^L$ of the graph (2.4) on the upper side
 of $\Sigma^L_g$ in such a way that the graph
is homotopic to the core $c^L$.  Analogously,
 we embed  the second copy $G^R$ of the graph
 (2.4) on the back side of $\Sigma^R_g$.

%introduce a copy $G^L$, resp. $G^R$, of the graph (4) on
%$\Sigma^L_g$, resp. $\Sigma^R_g$, in such a way that it is homotopic
%in the respective tubular neighborhood to its core $c^L$, resp. $c^R$,
%and lies on the upper side of $\Sigma^L_g$, resp. back side of $\Sigma^R_g$.

Finally, we make $G^L$ lefthanded and $G^R$ righthanded, i.e. we
introduce  meridians on each of the tubes corresponding to the lines
of $c^L$, resp. $c^R$, which undercross, resp. overcross, the lines
of $G^L$ on $\Sigma^L_g$, resp. $G^R$ on $\Sigma^R_g$. We then define
\be \label{K} K_{\underline{e},\underline{f}} =
\sum_{\underline{x},\underline{y}}
\prod^{3g-3}_{i=1}\frac{\omega^2_{x_i}}{\omega^2}
\frac{\omega^2_{y_i}}{\omega^2}\; Z(\tilde{M_g}, G^L_{\underline{e}}
\cup m^L_{\underline{x}} \cup G^R_{\underline{f}} \cup
m^R_{\underline{y}} \cup G^g)\ee
where $\underline{e}$, resp. $\underline{f}$, is a colouring of $G^L$,
resp. $G^R$, the product is over meridians and the sum is over
colourings $\underline{x}$ and
$\underline{y}$ of the meridians $m^L$ and $m^R$, on $\Sigma^L_g$ and
$\Sigma^R_g$, respectively, and $G^g$ is some canonical graph on $\Sigma_g$
without over- or undercrossings.

We denote by $V^L_g$, resp. $V^R_g$, the vector space associated to
$G^L$, resp. $G^R$, regarded as embedded into $\Sigma^L_g$,
resp. $\Sigma^R_g$, i.e.
\be V^L_g = \oplus_{\underline{e}} V^L_g(\underline{e}),\ee
where $V^L_g(\underline{e})$ is
 the tensor product of vector spaces associated to the coloured
3-vertices of $G^L$ taking into account the orientation of
$\Sigma^L_g$ and similarly for $G^R$.  Then
\be dim V^L_g = dim V^R_g = tr({\vec{N}}^2)^{(g-1)}\ee
by a simple counting, and hence
\be dim (V^L_g\otimes V^R_g) = dim V_{\Sigma_g}.\ee
Moreover, with the chosen orientation convention we have (see [BD])
$K_{\underline{e},\underline{f}} \in V^L_g(\underline{e})^\ast \otimes
V^R_g(\underline{f})^\ast \otimes
V_{\Sigma_g}$ and hence (\ref{K}) defines an operator
$$ K_{\underline{e},\underline{f}}: V^L_g(\underline{e}) \otimes
V^R_g(\underline{f}) \rightarrow V_{\Sigma_g}$$
in an obvious way.
%The basis in $V_\Sigma$ was actually constructed in [BD] by introducing
%two tubes (a right- and a lefthanded one) inside a handlebody
%$M_{\Sigma^g}$. The core of the tubes  is depicted below
%\be\label{kanon}
%\mbox{\epsfysize=2.5cm\epsffile{kanonnew.ps}}
%\ee
%and they are embedded in $M_{\Sigma^g}$ in such a way that
%the circles
%are non-contractible in $M_{\Sigma^g}$ and the tubes do not intersect
%each other.
%A tube called lefthanded (resp. righthanded) if it contains on its
%boundary a graph $G$ homotopic to the core of the tube and each line
%of the graph $G$ overcrosses (resp. undercrosses) exactly once a meridian
%$m$ on the tube. A state sum of a manifold $M$ with a lefthanded (or
%righthanded) tube $T$ inside is defined as follows:
%\be \label{nicht}
%Z(M(T), G_{\underline{a}})=\sum_{i,j, ...}
%\frac{\omega^2_i}{\omega^2}
%\frac{\omega^2_j}{\omega^2} ... Z(M(T), G_{\underline{a}} \cup m_i\cup
%m^\prime_j\cup ...)\ee
%where the sum is over colours of the meridians $m$, $m^\prime$, ...  .
We intend to show that the direct sum over
$\underline{e},\underline{f}$ of these operators yields an isomorphism
between $V^L_g \otimes V^R_g$ and $V_{\Sigma_g}$. This was proven for the
case $g=1$ in [BD]. In the general case it is a consequence of Lemma 1
below in which, however,
we have found it convenient first to rewrite $K_{\underline {e},
\underline{f}}$, up to a factor $\omega^{2(-g+1)}$, as
\be K_{\underline{e}, \underline{f}} = \sum_{\underline{x}}
\prod^{3g-3}_{i=1} \frac
{\omega^2_{x_i}}{\omega^2}\; Z(M^\prime_g,
G_{\underline{e},\underline{f}} \cup m_{\underline{x}} \cup G^g),\ee
where $M^\prime_g$ is the manifold with boundary components $\Sigma_g$
and ${\Sigma^{\prime}_g}^\ast$ obtained by removing one tubular neighborhood
instead of two as above and where $G_{\underline{e},\underline{f}}$ is
the coloured graph on $\Sigma^\prime_g$ indicated on the figure below
together with a system $m$ of meridians (of which there are $3g-3$
for $g\ge 1$, and 1 for $g=1$), and $G^g$ is as above.
%Observing that the  identifying of the boundaries of two tubes
%will change the state sum by a factor only,
%we will take  here the following canonical graph
%$G^g_{\underline{e}\,\underline{f}}$ as providing a basis in
%$V_{\Sigma_g}$:

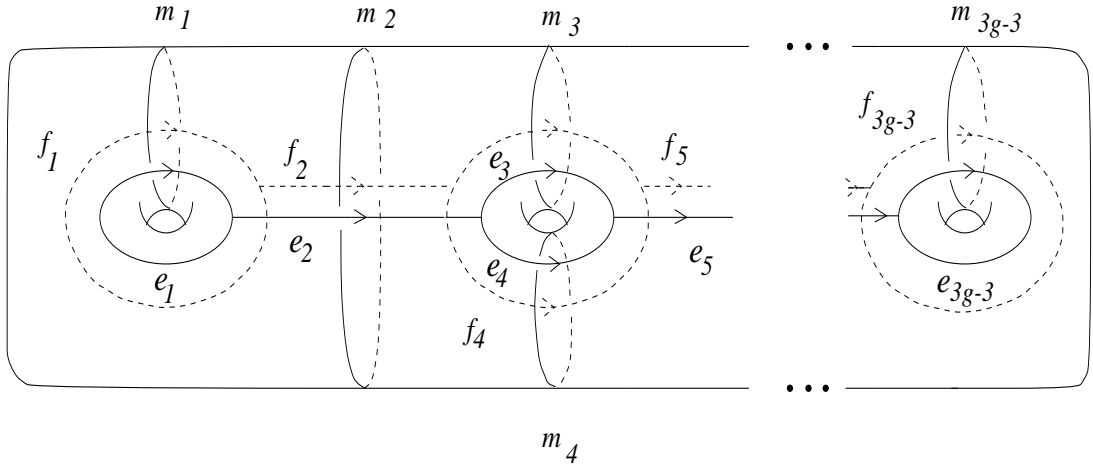
\be \label{basis}
%\mbox{\epsfysize=7cm\epsffile{basis.ps}}
\centerline{\hbox{
\psfig{figure=basis.ps,height=6.5cm,width=14.4cm}}}
\ee

The equivalence of (2.5) and (2.9) follows by merging $\Sigma^L_g$ and
$\Sigma^R_g$ as in the proof of Lemma 4.4 {\it ii)} in [BD]; see also
the proof of Lemma 1 below, where the same technique is used. We shall
henceforth take (2.9)
as the definition of $K_{\underline{e},\underline{f}}$.

%Letting ${\bar G}_{\underline{e},\underline{f}}$ denote the graph
%obtained from $G_{\underline{e},\underline{f}}$ by applying an
%orientation preserving element
%$S$ of the mapping class group of $\Sigma^\prime_g$ which interchanges the
%a- and b-cycles in a canonical homology basis for $\Sigma^\prime_g$ we
%define similarly

%We shall view

We now introduce an operator
$$ L_{\underline
{e},\underline{f}}: V_{\Sigma_g} \rightarrow
V^L_g(\underline{e})\otimes V^R_g(\underline{f})\; \subseteq
V^L_g\otimes V^R_g \, $$
as a mirror image of $K_{\underline{e},\underline{f}}$
w.r.t. a plane parallel
to the  $z$-axis and not intersecting the handlebody
$M_g$. More precisely,
\be L_{\underline{e},\underline{f}} = \sum_{\underline{x}}
\prod^{3g-3}_{i=1}\frac{\omega^2_{x_i}}{\omega^2}\;
Z(M^{\prime\prime}_g, {\bar G}_{\underline{e},\underline{f}}\cup
{ m}_{\underline{x}} \cup {\bar{G}}^g),\ee
where $M^{\prime\prime}_g$ is the mirror image of $M^{\prime}_g$ and $\partial
M^{\prime\prime}_g= \Sigma^\ast_g\cup \Sigma^{\prime\prime}_g$.
The graphs
 $  {\bar G}_{\underline{e},\underline{f}}\in  {\Sigma^{\prime\prime}_g}^\ast$
and
$ {\bar{ G}}_g\in \Sigma^\ast_g$ are the mirror images of
$   G_{\underline{e},\underline{f}}\in \Sigma^\prime_g$ and $G_g\in \Sigma_g$
 respectively.

Gluing $(M^\prime_g, G_{\underline{e},\underline{f}}\cup
m_{\underline{x}})$ and  $(M^{\prime\prime}_g,
{\bar G}_{\underline{e}^\prime,\underline{f}^\prime}\cup {
m}_{\underline{y}})$ along $\Sigma_g$ we
obtain $(N_g, G_{\underline{e},\underline{f}}
\cup m_{\underline{x}}, {\tilde{
G}}_{\underline{e}^\prime,\underline{f}^\prime}\cup {
m}_{\underline{y}})$ where $N_g$ is
diffeomorphic to $\Sigma_g\times [0,1]$ with boundary
$\Sigma^{\prime\prime}_g \cup {\Sigma^{\prime}_g}^\ast$.
%The only difference of
The graph
${\tilde{G}}_{\underline{e}^\prime,\underline{f}^\prime}
\cup m_{\underline{y}}
\in \Sigma^{\prime\prime}_g$ can be obtained from the standard
graph $G_{\underline{e},\underline{f}}
\cup m_{\underline{x}}
\in \Sigma^{\prime}_g$ depicted in (2.10)
by changing the colourings $\underline{e}\rightarrow \underline{e}^\prime$, $
\underline{f}\rightarrow \underline{f}^\prime$, $ \underline{x}\rightarrow
\underline{y}$ and replacing all overcrossings
 by undercrossings and vice versa.

%is that the $\underline{e}^\prime$-coloured lines
%go under the meridians  $m$
%and the $\underline{f}^\prime$-coloured lines over them.

Eq. (2.1) implies that
\be\label{sthand}
L_{\underline{e}^\prime,\underline{f}^\prime}K_{\underline{e},\underline{f}}
= \sum_{\underline{x},\underline{y}}
\prod_i\frac{\omega^2_{x_i}}{\omega^2}\frac{\omega^2_{y_i}}{\omega^2}\;
Z(N_g, G_{\underline{e},\underline{f}}\cup
m_{\underline{x}}\cup {\tilde{
G}}_{\underline{e}^\prime,\underline{f}^\prime}\cup
m_{\underline{y}}) \, .\ee

We are now in position to state the announced lemma.

\begin{lem}  The operator
$L_{\underline{e}^\prime,\underline{f}^\prime}K_{\underline{e},\underline{f}}:
V^L_g(\underline{e}) \otimes V^R_g(\underline{f}) \rightarrow
V^L_g(\underline{e}^\prime)\otimes V^R_g(\underline{f}^\prime)$
satisfies
\be\label{id}
\omega^{2g-2}
\omega_{\underline{e}}\omega_{\underline{f}}\omega_{\underline{e}^\prime}\omega_{\underline{f}^\prime}L_{\underline{e}^\prime,\underline{f}^\prime}K_{\underline{e},\underline{f}}
=
\delta_{\underline{e},\underline{e}^\prime}\delta_{\underline{f},\underline{f}^\prime}\;
{\eins}_{V^L_g(\underline{e})\otimes V^R_g(\underline{f})},\ee
where we have introduced the notation $\omega_{\underline{e}} =
\prod^{3g-3}_{i=1} \omega_{e_i}$ and
$\delta_{\underline{e},\underline{e}^\prime} =
\prod^{3g-3}_{i=1}\delta_{e_i,e^\prime_i}$.

\end{lem}
\vspace*{0.2cm}

{\bf Proof:} The idea of the argument is the following. By introducing
tubes between $\Sigma^{\prime}_g$ and $\Sigma^{\prime\prime}_g$ we
  step by step lift the  lines of
$G^g_{\underline{e} \underline{f}}\in \Sigma^\prime_g$
 on $\Sigma^{\prime \prime}_g$ and  cut
the  handles traversed by these lines.
 Applying the technique developed in [BD] and [KS]
we will arrive on (2.13).

%reduce the
%state sum in (\ref{sthand}) to the state sum for a handlebody
%$\tilde {N_g}$ with the graph $G^g_{\underline{e}\, \underline{f}}\cup
%G^g_{\underline{e}^\prime\, \underline{f}^\prime}$ and a set of
%meridians on
%its boundary $\tilde{\Sigma}_g$. By
%cutting the handles of $\tilde {N_g}$ we shall arrive at (\ref{id}).

 Due to Lemma 3.3 in [BD] introduction of a tube
with an $a$-coloured meridian (which is {\it not} normalized by
$\omega^{-2}$) does not change the state sum. Pictorially this looks as
follows:
\begin{center}
\mbox{\epsfysize=5cm\epsffile{tube1.ps}}
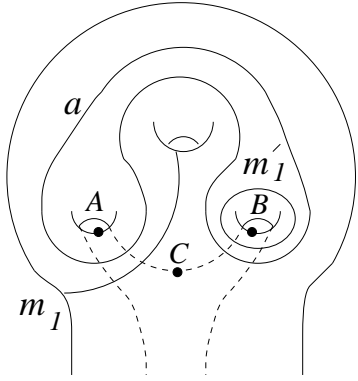
\\Fig.2 {\it A part of the manifold $N_g$ where the
boundary component $\Sigma^\prime_g$ of the tube is connected to
$\Sigma^{\prime\prime}_g$ by a  tube
with an $a$-coloured meridian on it}
\end{center}
where we do not draw the $\underline{e}$-, $\underline{f}$- and
$\underline{e}^\prime$-, $\underline{f}^\prime$-coloured lines.
Applying Lemma 4.2 {\it ii)} in [BD] (or the Wigner-Eckart type
relation (A.15) in [KS]) to the meridians
$m_1$, $m^\prime_1$ and  $a$  we can change the
graph so that
 the handle $(ABC)\times I$ will be traversed by a single line only.
 According to Remark 3.6 in [BD] the colour of this line
can be set to zero and the handle  cut. This yields a manifold
$N^\prime_g$  as depicted on Fig.3.
\begin{center}
\mbox{\epsfysize=7cm\epsffile{tube2.ps}}
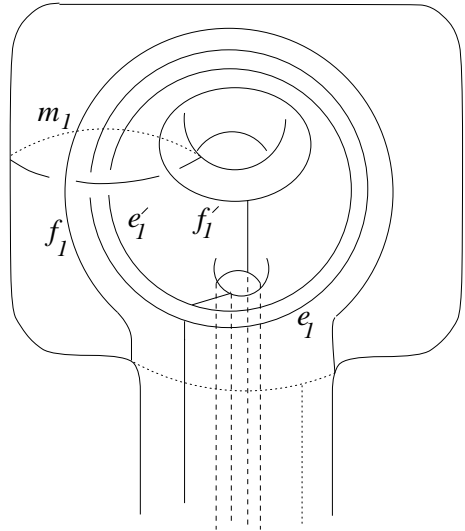
\\Fig.3 {\it A part of the manifold $N^\prime_g$ with associated graph on it}
\end{center}

Using  lemma 4.2 {\it ii)} in [BD] once more (see also example 5.8 {\it iii)}
 in  [KS]) one can
cut the handle traversed by $e^\prime_1$- , $e_1$- ,
$f^\prime_1$- and $f_1$-coloured lines. After that the state sum of the
resulting $(g-1)$-cylinder
becomes multiplied by $\omega^{-2}_{e_1}\omega^{-2}_{f_1} \delta_{e^\prime_1
e_1}\delta_{f^\prime_1 f_1}$.

Continuing this procedure analogously
we obtain
%L_{\underline{e^\prime}\,\underline{f^\prime}}K_{\underline{e}\,\underline{f}}\,=
%\frac{1}{\omega^{2(g-1)}} \sum_{\underline{x}} \prod_i
%\frac{\omega^2_{x_i}}{\omega^2}
%Z(\tilde {N_g}, G^g_{\underline{e}\,\underline{f}}
%\cup G^{g}_{\underline{e}^\prime\,\underline{f}^\prime}\cup
%%m_{\underline{x}})\, .$$
%Removing now the handles of $\tilde {N_g}$ using (4.4) and (4.9-11) in
%[BD] (see also example
%5.8 {\it iii)} in [KS]) we obtain
the desired result:
$$L_{\underline{e}^\prime\,\underline{f}^\prime}K_{\underline{e}\,\underline{f}}=
\omega^{-2g+2}
\;\delta_{\underline{e}\,\underline{e}^\prime}
\delta_{\underline{f}\,\underline{f}^\prime}\; (\omega^2_{\underline{e}}\;
\omega^2_{\underline{f}})^{-1}\; {\eins}_{V^L_g(\underline{e})\otimes
V^R_g(\underline{f})}\,  .$$

$\hfill\Box$
\vspace*{0.2cm}

%2) Show that
%the matrix $S\otimes S^\ast$ with the elements
%\be\label{smatr1}
%\mbox{\epsfysize=4cm\epsffile{smatr1.ps}}
%\ee
%and $S^\ast_{\underline{e}\, \underline{e^\prime}} $
%obtained from the def.(\ref{smatr1}) by interchanging all over- and
%undercrossings, defines a permutation of all
%contractible and non-contractible cycles in the basis (\ref{basis})
%on $V_\Sigma$. See also that the matrix $S$ is unitary, i.e.\
%$$\sum_{\underline{e^\prime}}\,S_{\underline{e}\, \underline{e^\prime}} \;
%S^\ast_{\underline{e^\prime}\, \underline{e^{\prime\prime}}}
%=\delta_{\underline{e}\,\underline{e^{\prime\prime}}}\, .$$.

%Help: Calculate that
%$$\omega^{2g-2}\,\omega_{\underline{e}}\,\omega_{\underline{e^\prime}}\,
%\omega_{\underline{f}}\, \omega_{\underline{f^\prime}}\,
%Z(\Sigma_g\times I,
%G^g_{\underline{e}\,\underline{f}}\cup
%G^g_{\underline{e^\prime}\,\underline{f^\prime}})=
%S_{\underline{e}\, \underline{e^\prime}}
%\; S^\ast_{\underline{f}\, \underline{f^\prime}} $$
%where $G^g_{\underline{e}\,\underline{f}}\subset \Sigma\times \{0\}$
%identifyed with $\Sigma$ and
%$G^g_{\underline{e^\prime}\,\underline{f^\prime}}\subset \Sigma\times
%\{1\}$ diffeomorphic to $\Sigma^\ast$.

Defining the operators $K: V^L_g \otimes V^R_g \rightarrow
V_{\Sigma_g}$  and $L: V_{\Sigma_g} \rightarrow V^L_g \otimes V^R_g$
by
$$ K = {\omega}^{g-1} \oplus_{\underline{e},\underline{f}}
\omega_{\underline{e}} \omega_{\underline{f}}
K_{\underline{e},\underline{f}}\;, \;  \; L = {\omega}^{g-1}
\oplus_{\underline{e},\underline{f}} \omega_{\underline{e}}
\omega_{\underline{f}} L_{\underline{e},\underline{f}}$$
it follows from (\ref {id}) that $LK = {\eins}_{V^L_g \otimes V^R_g}$ and
hence by (2.8) $K$ and $L$ are isomorphisms and
\be\label{iso} L = K^{-1}\;  .\ee

Although we shall strictly speaking not use them in the following let us
introduce the left- and righthanded counterparts
$K^L_{\underline{e}}$ and $K^R_{\underline{f}}$ of
$K_{\underline{e},\underline{f}}$ by replacing in eq. (2.9) the graph
$G_{\underline{e},\underline{f}}$ by its left- and righthanded parts
$G^L_{\underline{e}}$ and $G^R_{\underline{f}}$, respectively, and
similarly $L^L_{\underline{e}}$ and $L^R_{\underline{f}}$ by replacing
 ${\bar G}_{\underline{e},\underline{f}}$ in eq. (2.11) by
${\bar G}^L_{\underline{e}}$ and ${\bar G}^R_{\underline{f}}$, respectively.
The
proof of Lemma 1 then yields
$${\omega}^{2g-2}\;  \omega_{\underline{e}}\;
\omega_{\underline{e}^\prime}\;
L^L_{\underline{e}} K^L_{\underline{e}^\prime} =
{\delta}_{\underline{e},\underline{e}^\prime}
{\eins}_{V^L_g(\underline{e})}$$
and
$${\omega}^{2g-2}\;  {\omega}_{\underline{f}}\;
{\omega}_{\underline{f}^\prime}\;  L^R_{\underline{f}} K^R_{f^\prime} =
{\delta}_{\underline{f},\underline{f}^\prime}
{\eins}_{V^R_g(\underline{f})}$$
and consequently
$$ L^L K^L = {\eins}_{V^L_g}\; , \;  \;  L^R K^R = {\eins}_{V^R_g}\;
,$$
where $K^L:V^L_g \rightarrow V_{\Sigma_g}$ and $L^L: V_{\Sigma_g}
\rightarrow V^L_g$ are defined by
\be K^L = {\omega}^{g-1} \oplus_{\underline{e}}
{\omega}_{\underline{e}} K^L_{\underline{e}}\;  ,\;  \; L^L =
{\omega}^{g-1} \oplus_{\underline{e}} {\omega}_{\underline{e}}
L^L_{\underline{e}},\ee
and similarly for $K^R: V^R_g \rightarrow V_{\Sigma_g}$ and $L^R:
V_{\Sigma_g} \rightarrow V^R_g$ .

\section{Factorization of state sums}
\label{fac}  For each
genus $g \geq 0$ we fix once and for all manifolds $M^\prime_g$ and
$M^{\prime\prime}_g$ as defined in Section 2 with
$\partial{M^\prime_g}
= \Sigma_g \cup {\Sigma^{\prime}_g}^\ast$ and $\partial{M^{\prime\prime}_g} =
\Sigma^\ast_g \cup {\Sigma^{\prime\prime}_g}$ , where $\Sigma_g$,
$\Sigma^\prime_g$ and $\Sigma^{\prime\prime}_g$ are fixed oriented
surfaces of genus $g$ and where fixed graphs $G^g_{\underline{e} ,
\underline{f}}$ and  $\bar{
G}^g_{\underline{e},\underline{f}}$ are embedded in
$\Sigma^\prime_g$ and ${\Sigma^{\prime\prime}_g}^\ast$
respectively, together with the
associated sets of meridians.
%$m^g$, resp. $\bar{m}^g$.
We have here
made the dependence of the graphs and meridians on the genus explicit,
and will do so likewise for the associated operators
$K_{\underline{e},\underline{f}},
L_{\underline{e},\underline{f}}$ etc.

By a parametrized surface of genus $g$ we  mean a pair $(\Sigma,\phi)$,
where $\Sigma$ is a compact, connected, oriented surface of genus $g$ and
$\phi: \Sigma
\rightarrow \Sigma_g$ is a diffeomorphism. We call $\phi$ a
parametrization of $\Sigma$
%In case $\phi: \Sigma \rightarrow
%\Sigma_g$ is orientation reversing, i.e. $\phi:\Sigma \rightarrow
%\Sigma^\ast_g$ is orientation preserving,
and set
$$ \tilde{V}_\Sigma(\phi) = V^L_g \otimes V^R_g \, . $$
Let us
consider a 3-dimensional cobordism $M$
whose  boundary
$\partial M = \Sigma^\ast_1 \cup\Sigma_2$
consists of two  compact, connected, oriented
 surfaces of genus $g_1$ and $g_2$, respectively, which are
parametrized by $\phi_1$ and $\phi_2$.
  An operator
$\tilde{Z}(M):\tilde{V}_{\Sigma_1 }(\phi_1)
\to\tilde{V}_{\Sigma_2 }(\phi_2)$
can be defined as follows:
$$\tilde{Z}(M)= L(\phi_2)Z(M)K(\phi_1)\, ,$$
where
$$K(\phi_1) = U(\phi_1) K^{g_1}, \;\;\;\;\;
 L(\phi_2)= L^{g_2} U(\phi_2)$$
and $U(\phi  ):V_\Sigma \to V_{\Sigma_g}$ satisfying (1.2).

%\equiv and denote by $W(\phi): V_\Sigma \rightarrow \tilde{V}_\Sigma(\phi)$
%the isomorphism given by
%\be \label{w} W(\phi) = (K^g)^t U(\phi), \ee On the other hand,
%if $\phi:\Sigma \rightarrow \Sigma_g$ is orientation preserving we set
%$$\tilde{V}_\Sigma(\phi) = V^L_g \otimes V^R_g$$ and define the
%isomorphism $W(\phi): V_\Sigma \rightarrow \tilde{V}_\Sigma(\phi)$ by
%

%More generally, let us consider the cobordism $M$
%with boundary components $ \Sigma^1_{g_1}, ... ,
% \Sigma^n_{g_n}$ and parametrization $\phi_i$
%(\phi_i)$$ and $$ W(\phi_1,...,\phi_n) = \otimes^n_{i=1} W(\phi_i),$$
%\Sigma_g$ denotes the corresponding parametrization of $\Sigma^\ast$
%then
%\be \label{transw} W(\phi^\ast) = (W(\phi)^t)^{-1} \ee
%as a consequence of
%(\ref{transu}), (\ref{w}), (\ref{w1}) and (\ref{iso}).

More generally,
given a compact, oriented cobordism $M$ with boundary
components
$ {\Sigma^1_{g_1}}^\ast , ... ,{\Sigma^m_{g_m}}^\ast ,
\Sigma^{m+1}_{g_{m+1}}, ..., \Sigma^n_{g_n}$
and parametrization $\phi_i$ of  $\Sigma^i_{g_i}$
we set
\be\label{npart}
\tilde{Z}(M)= L(\phi_{m+1}, ...,\phi_{n}) Z(M) K(\phi_1, ...,\phi_n)\ee
where
 $$ K(\phi_1,...,\phi_k) = \otimes^k_{i=1} K(\phi_i)$$
and $ L(\phi_1,...,\phi_k)$
% = \otimes^k_{i=1} L(\phi_i)\, .$$
is defined analogously.
%$\partial{M} = \Sigma$ as given above we define
%\be \label{npart} \tilde{Z}(M) = W(\phi_1,...,\phi_n)(Z(M)). \ee

Equivalently, (\ref{npart}) can be expressed as follows. Let $\bar M$
denote the
manifold obtained by gluing $M^\prime_{g_i}$ onto $M$
along $\phi_i$ for $1<i<m$, and gluing
$M^{\prime\prime}_{g_i}$ onto $M$ along $\phi_i$ in case
$i>m$.
%Let us for convenience order the boundary components such that
%$\phi_1,...,\phi_m$ are orientation reversing, wheras $\phi_{m+1},...,\phi_n$
Then, clearly, $\bar M$ is diffeomorphic to
$M$ and has boundary components
$(\Sigma^{\prime}_{g_1})^\ast ,...,(\Sigma^{\prime}_{g_m})^\ast ,$
$ \Sigma^{\prime\prime}_{g_{m+1}},..., \Sigma^{\prime\prime}_{g_n}$ with
embedded graphs
$G^{g_1}_{\underline{e}^1,\underline{f}^1},...,G^{g_m}_{\underline{e}^m,\underline{f}^m},\-
{\bar G}^{g_{m+1}}_{\underline{e}^{m+1},\underline{f}^{m+1}},...,
{\bar G}^{g_n}_{\underline{e}^n,\underline{f}^n}$, respectively.
With the notation $\tilde{e} = (\underline{e}^1,...,\underline{e}^n)$
and  $$\omega_{\tilde{e}} = \prod^n_{i=1} \omega_{\underline{e_i}}$$ we
then have
\be\  \tilde{Z}(M) =
\oplus_{\tilde{e},\tilde{f}}\tilde{Z}_{\tilde{e},\tilde{f}}(M), \ee
where the coloured state sum $\tilde{Z}_{\tilde{e},\tilde{f}}(M)$ is
defined by
\begin{eqnarray} \label{ztilde}
\tilde{Z}_{\tilde{e},\tilde{f}}(M) =
\omega^{g_1 +...+ g_n -n}
\omega_{\tilde{e}}\;  \omega_{\tilde{f}} \; \sum_{\tilde{x}}
\prod_{i,j} \frac{\omega^2_{x^j_i}}{\omega^2} Z(\bar M,
{\cal{G}}_{\tilde{e},\tilde{f}} \cup { \cal M}_{\tilde{x}})
\end{eqnarray}
where $$ {\cal{G}}_{\tilde{e},\tilde{f}} =
G^{g_1}_{\underline{e}^1,\underline{f}^1} \cup
... \cup {\bar{G}}^{g_n}_{\underline{e}^n, \underline{f}^n} $$
and $$ {\cal M}_{\tilde{x}} = m^1_{\underline{x}^1} \cup ... \cup
m^n_{\underline{x}^n} \, .$$
%and the full state sum
%\be \tilde{Z}(M) = \oplus_{\tilde{e},\tilde{f}}
%\tilde{Z}_{\tilde{e},\tilde{f}}(M)\;  . \ee
%Equivalently, if we consider $M$ as a cobordism and interpret the
%state sum as an operator $Z(M): V^\ast_{\Sigma^1_{g_1} \cup ...\cup
%\Sigma^m_{g_m}} \rightarrow V_{\Sigma^{m+1}_{g_{m+1}}\cup ...\cup
%%\Sigma^n_{g_%n}}$  then
%$Z_{\tilde{e},\tilde{f}}(M): \otimes^m_{i=1}
%(V^L_{g_i}(\underline{e}^i)^\ast \otimes
%V^R_{g_i}(\underline{f}^i)^\ast) \rightarrow \otimes^n_{j=m+1}
%(V^L_{g_j}(\underline{e}^j) \otimes V^R_{g_j}(\underline{f}^j))$ is
%given by
%\be \label{til} \tilde{Z}_{\underline{e},\underline{f}}(M) =
%\omega^{g_1+...g_n -n} \omega_{\tilde{e}} \omega_{\tilde{f}}
%\otimes^n_{j=m+1} (L^{g_j}_{\underline{e}^j,\underline{f}^j}
%U(\phi_j)) Z(M) \otimes^m_{i=1} (U(\phi_i)^{-1}
%K^{g_i}_{\underline{e}^i,\underline{f}^i})\; , \ee
%where we note that $\phi_i: (\Sigma^i_{g_i})^\ast \rightarrow
%\Sigma_{g_i}\; , i=1,...,m, \;$  and $ \phi_j: \Sigma^j_{g_j}
%\rightarrow \Sigma_{g_j},\;  j=m+1,...,n,\; $  are all orientation
%preserving. Similarly, we have
%\be \label{till} \tilde{Z}(M) = \otimes^m_{j=1} (L^{g_j} U(\phi_j))
%Z(M) \otimes^n_{i=m+1} (U(\phi_i)^{-1} K^{g_i})\;.\ee

Finally, we define an isomorphism $ \tilde{U}(f):
\tilde{V}_\Sigma(\phi) \rightarrow
\tilde{V}_{\Sigma^\prime}(\phi^\prime)$ by
\be \label{utilde} \tilde{U}(f) = L(\phi^\prime) U(f) K(\phi),\ee
for any orientation preserving diffeomorphism $f:\Sigma \rightarrow
\Sigma^\prime$ between pa\-ra\-me\-trized surfaces $(\Sigma,\phi)$ and
$(\Sigma^\prime,\phi^\prime)$ of genus $g$. This definition is
extended in an obvious way to orientation preserving diffeomorphisms
between arbitrary compact, oriented surfaces in terms of tensor
products.

The objects $\tilde{V}, \tilde{U},\tilde{Z}$ define a TQFT on compact,
oriented 3-manifolds with parametrized boundary.
This can be easily verified using
the definition of these objects and eq. (2.14).
%except for the gluing property 3) of
%Section 1, which is a consequence of Lemma 1. Indeed, let $(\Sigma,
%\phi)$ and $(\Sigma^\prime,\phi^\prime)$ be two parametrized surfaces
%of genus $g$ and let $M_1$ and $M_2$ be compact, oriented 3-manifolds
%with parametrized boundaries containing $(\Sigma,\phi)$ and
%$({\Sigma^\prime}^\ast,{\phi^\prime}^\ast)$, respectively, and such that
%$\phi^\prime \phi^{-1}$ is orientation preserving, i.e. either $\phi$
%and $\phi^\prime$ are both orientation preserving or they are both
%orientation reversing. The parametrization isomorphisms $W_1$ and
%$W_2$ associated
%with $\partial{M_1}$ and $\partial{M_2}$, respectively, can then be
%written in factorized form as $W_1 = \bar{W}_1 \otimes
%W(\phi)$ and $\bar{W}_2 \otimes
%W({\phi^\prime}^\ast)$ such that, given an orientation preserving
%diffeomorphism
%$f: \Sigma \rightarrow \Sigma^\prime$, we have, in operator language,
%\begin{eqnarray*}
% \tilde{Z}(M_2) \tilde{U}(f) \tilde{Z}(M_1) & = & \bar{W}_2 Z(M_2)
%W({\phi^\prime}^\ast)^t \tilde{U}(f) W(\phi) Z(M_1) \bar{W}_1^t \\
%& = &  \bar{W}_2 Z(M_2) U(f) Z(M_1) \bar{W}_1^t \\
%& = & \bar{W}_2 Z(M) \bar{W}^t_1 = \tilde{Z}(M),
%\end{eqnarray*}
%where in the second step we have used
%the definition of $\tilde{U}(f)$ and (\ref{transw}), and $M$ denotes
%the manifold obtained by gluing $_1$ and $M_2$ along $f$.
%Note also that $\tilde{U}(f)$ satisfies the condition (\ref{transu})
%as a consequence of (\ref{transw}).
The TQFT based on $\tilde{V}, \tilde{U}$ and $ \tilde{Z}$
 is equivalent to the theory defined in the
previous section. The equivalence is given by the $K$ and $ L$-operators
(see [T] or [DJ]).

We are now ready to state and prove the main result of this paper.

\begin{satz}
Let $M$ be a compact, oriented 3-manifold. For any colouring
$(\tilde{e},\tilde{f})$ as defined above we have
\be \label{tens} \tilde{Z}_{\tilde{e},\tilde{f}}(M) \, =\,
\tau_{\tilde{e}}(M)\otimes
\tau_{\tilde{f }} (M^\ast)\ee
where the invariant $\tau_{\tilde{e}}$ is given by
eq. (\ref{end}) below and coincides with the invariant introduced in
[T] up to normalization.

\end{satz}

{\bf Proof:}

As remarked earlier, we can replace each tube in ${\bar{M}}$ defined
above with
graph $G^{g_i}_{\underline{e}^i,\underline{f}^i} \cup
m^i_{\underline{x}^i}$ by two tubes with graphs
$(G^{g_i}_{\underline{e}^i})^L \cup (m^i_{\underline{x}^i})^L$ and
$(G^{g_i}_{\underline{f}^i})^R \cup (m^i_{\underline{y}^i})^R$,
respectively, at the cost of a factor $\omega^{2(g_i-1)}$. Let us assume we
have done so for each $i = 1,...,n$
 and denote the resulting manifold also by ${\bar{M}}$.
As is well known, the closed
manifold obtained from ${\bar{M}}$ by filling  all $2n$ tubes
%that were removed by the definition of $M^\prime_{g_i}, i=1,...,n$,
has a representation
 by surgery on $S^3$ along a set of links $l_1,...,l_N$
which, of course, may be assumed not to intersect the filled  tubes.
%The manifold $\tilde{S}^3$ obtained
%by removing disjoint tubular neighborhoods $T_1,...,T_N$ of
%$L_1,...,L_N$ from $\bar{\bar{M}}$ and
Using Lemma 1 for the case $g = 1$ as in the proof of Theorem 5.2 in
[BD] one obtains
\begin{eqnarray}\label{arr} &   & Z({\bar{M}},
G^{g_1}_{\underline{e}^1,\underline{f}^1} \cup m^1_{\underline{x}^1}
\cup ... \cup \bar{G}^{g_n}_{\underline{e}^n,\underline{f}^n} \cup
m^n_{\underline{x}^n})= \nonumber \\ & = &
 \omega^{2(g_1+...+g_n-n-N)}
\sum_{\tilde{a},\tilde{z},\tilde{b},\tilde{z}^\prime}
\omega^2_{\tilde{a}} \; \omega^2_{\tilde{b}}
\; \frac{\omega_{\tilde{z}}^2}{\omega^{2N}}
 \; \frac{\omega_{\tilde{z}^\prime}^2}{\omega^{2N}} \nonumber \\
&   &  Z(\tilde{S}^3, {\cal L}^L_{\tilde{a}} \cup
({\cal{M}}^\prime_{\tilde{z}})^L \cup {\cal{L}}^R_{\tilde{b}} \cup
({\cal{M}}^\prime_{\tilde{z}^\prime})^R \cup
{\cal{G}}^L_{\tilde{e}} \cup {\cal M}^L_{\tilde{x}} \cup
{\cal G}^R_{\tilde{f}} \cup {\cal M}^R_{\tilde{y}})
\end{eqnarray}
where we have introduced the shorthand notation $${ \cal G}^L_{\tilde{e}} =
(G^{g_1}_{\underline{e}^1})^L \cup ... \cup
(\bar{G}^{g_n}_{\underline{e}^n})^L$$
and similarly for the righthanded part and the meridians. Furthermore,
$ \tilde{S}^3$ denotes the manifold obtained from ${\bar{M}}$ by
removing  two disjoint tubular neighborhoods $ T^L_i$ and
$T^R_i$ for each $i = 1,...,N$. We define
%The neighborhoods
$ T^L_i$ and $T^R_i$
 by splitting a tubular neighborhood of $l_i$ into two
nearby ones as was done previously for the graphs
$G^{g_1},...,G^{g_n}$.
Finally, ${\cal L}^L = {L_1}^L \cup ... \cup {L_N}^L$
(together with  associated
meridians ${\cal M}^L = m^L_1 \cup ... \cup m^L_N$)
is a collection
of lefthanded graphs
on the boundary components $\partial T^L_1,..., \partial T^L_N$
of $\tilde{S}^3$,
where the graphs are determined by the surgery prescription, and
similarly for $\cal{L}^R$ and $\cal{M}^R$.

Next we recall from [BD] (see also [KS]) that two tubes with left- and
righthanded lines, respectively, have trivial braiding, i.e. they may
be deformed through each other. Using this and the fact that
$\tilde{S}^3$ is a 3-sphere with a collection of $2(n+N)$ tubes
removed, together with the
factorization property of $Z(M,G)$ w.r.t. connected sums (see Lemma 3.2
in [BD]), we obtain by substituting  (\ref{arr}) into
(\ref{ztilde}) that
\be \label{fin} \tilde{ Z}_{\tilde{e},\tilde{f}}(M)
 = \omega^{2(g_1+...g_n-n-N+1)} \sum_{\tilde{a},\tilde{b}} { \cal Z}(S^3,
{\cal L}^L_{\tilde{a}} \cup {\cal G}^L_{\tilde{e}}) \otimes {\cal Z}(S^3,
{\cal L}^R_{\tilde{b}} \cup {\cal G}^R_{\tilde{f}}),
\ee
where we have introduced
\begin{eqnarray*} {\cal Z}(S^3,{\cal L}^L_{\tilde{a}} \cup
{\cal G}^L_{\tilde{e}}) &=& \omega^{g_m + ... + g_n -(n-m)}
\omega_{\tilde{e}}
 \omega^2_{\tilde{a}} \sum_{\tilde{z}, \tilde{x}}
\frac{\omega^2_{\tilde{z}}}{\omega^{2N}} \; \prod_{i,j}
 \frac{\omega^2_{x^j_i}}{\omega^2}\\ & &
Z((\tilde{S}^3)^L,{ \cal L}^L_{\tilde{a}} \cup
({\cal M}^\prime_{\tilde{z}})^L \cup {\cal G}^L_{\tilde{e}} \cup
{\cal M}^L_{\tilde{x}})
\end{eqnarray*}
and
\begin{eqnarray*}  {\cal Z}(S^3,{\cal L}^R_{\tilde{b}}
\cup
{\cal G}^R_{\tilde{f}}) &=& \omega^{g_1 + ... + g_m -m} \omega_{\tilde{f}}
 \omega^2_{\tilde{b}} \sum_{\tilde{z^\prime}, \tilde{y}}
\frac{\omega^2_{\tilde{z^\prime}}}{\omega^{2N}} \; \prod_{i,j}
 \frac{\omega^2_{y^j_i}}{\omega^2}\\ & &
Z((\tilde{S}^3)^R,{ \cal L}^R_{\tilde{b}} \cup
({\cal M}^\prime_{\tilde{z^\prime}})^R \cup {\cal G}^R_{\tilde{f}} \cup
{\cal M}^R_{\tilde{y}})
\end{eqnarray*}
where $(\tilde{S}^3)^L$ is defined in analogy with $\tilde{S}^3$
except that only tubes with lefthanded graphs or links are removed
from $S^3$ and $(\tilde{S}^3)^R$ is defined similarly.

Finally, setting $$\Delta_L = \sum_{c \in {\cal I}} q^2_c \omega^4_c ,$$
 we define
\be \label{end}  \tau_{\tilde{e}}(M) = \omega^{g_1+...+g_n-n-N+1} (\Delta_L
\omega^{-1})^{\sigma({\cal L})} \sum_{\tilde{a}}
{\cal Z}((\tilde{S}^3)^L, {\cal L}^L_{\tilde{a}} \cup
{\cal G}^L_{\tilde{e}}),
\ee
 where
 $\sigma(\cal{L})$ is the signature of a certain 4-manifold whose
boundary is ${\bar{M}}$ with tubes filled in. Similarly, the
righthanded counterpart $\tau^R_{\tilde{f}}$ is defined with
 $\Delta_R$ given by the same formula as $\Delta_L$ except that $q_c$ should
be replaced by $q^{-1}_c$. Then $$ \Delta_L \Delta_R = \omega^2$$ (see
[T]) and hence (\ref{fin}) can be rewritten as  $$
Z_{\tilde{e},\tilde{f}}(M) = \tau_{\tilde{e}}(M)
\otimes\tau^R_{\tilde{f}}(M).$$
By arguments identical to those in [BD] one shows that $$\tau^R_{\tilde{f}}(M)
= \tau_{\tilde{f}}(M^\ast)$$
 thus proving
(\ref{tens}). Likewise the argument that $\tau_{\tilde{e}}(M)$ equals
the ribbon graph invariant introduced in [T] follows as in [BD] by
projecting the tubes in $(\tilde{S}^3)^L$ with graphs and links onto a plane.
 $\hfill\Box$
\vspace*{0.2cm}

\section{Concluding remarks}

The proof of Theorem 2 can be extended in a straightforward manner to
the case  where
punctures are introduced on the boundary components of {M}. We shall,
however, not elaborate on that case here (see also [T]).

It should be mentioned that the equivalence of the TQFT defined in
section 2 and the one defined in terms of $\tilde{V}, \tilde{U},
\tilde{Z}$ follows from the equality of the corresponding state sums
of closed manifolds, shown in [BD] and [T], once it is known that the two
theories are
non-degenerate (see e.g. [T]). The method of this paper gives the
equivalence explicitly and at the same time prepares the ground for
the proof of (\ref{tens}).

\vspace*{1cm}

{\Large \bf Acknowledgements}
\vspace*{1cm}

One of us (A.B.) would like to thank  Vladimir Turaev for
valuable suggestions.

%One of us (A.B) would like to thank  Vladimir Turaev for
%suggesting the problem.

\end{document}